\documentclass[reqno,11pt,twoside,english]{amsart}
\usepackage[T1]{fontenc}
\usepackage[utf8]{inputenc}
\usepackage{lmodern,microtype}
\usepackage[english]{babel}
\usepackage{amsmath,amssymb,amsthm,mathrsfs}
\usepackage[final,colorlinks=true,linkcolor=blue,citecolor=blue,urlcolor=blue]{hyperref}

\allowdisplaybreaks[4]
\numberwithin{equation}{section}
\newcommand{\R}{\mathbb R}
\newcommand{\HH}{\mathcal H^1}
\newcommand{\hh}{h^1}

\newtheorem{theorem}{Theorem}[section]
\newtheorem{lemma}[theorem]{Lemma}
\newtheorem{proposition}[theorem]{Proposition}
\newtheorem{corollary}[theorem]{Corollary}
\theoremstyle{definition}

\newtheorem{remark}[theorem]{Remark}

\begin{document}

\title[Rivi\`ere's System with a Local Hardy Drift]{On $W^{2,1}$-Regularity for
Rivi\`ere's System with a Local Hardy Drift}

\author[Y.-T. Kang]{Yu-Ting Kang}
\address[Y.-T. Kang]{Research Center for Mathematics and Interdisciplinary Sciences, Shandong University 266237, Qingdao, P. R. China and Frontiers Science Center for Nonlinear Expectations, Ministry of Education, P. R. China}
\email{yuting.kang@mail.sdu.edu.cn}

\author[W.-J. Qi]{Wen-Juan Qi}
\address[W.-J. Qi]{The Department of Mathematical Sciences, Tsinghua University,
100084, Beijing, P. R. China}
\email{qiwenjuan2026@163.com}

\author[Z.-M. Sun]{Zhao-Min Sun}
\address[Z.-M. Sun]{Research Center for Mathematics and Interdisciplinary Sciences, Shandong University 266237, Qingdao, P. R. China and Frontiers Science Center for Nonlinear Expectations, Ministry of Education, P. R. China}
\email{sunzhaomin@mail.sdu.edu.cn}

\thanks{Corresponding author: Zhao-Min Sun.}

\begin{abstract}
In this paper, we study the endpoint regularity of weak solutions to the inhomogeneous Rivière system
$$
-\Delta u=\Omega\cdot\nabla u+f
\qquad\text{in }B_1\subset\mathbb R^2,
$$
with an antisymmetric potential $\Omega\in L^2$ and a  local Hardy drift term $f\in h^1$. We construct a continuous solution in $W_0^{1,2}$, which fails to belong to $W^{2,1}$ in every neighborhood of the origin, even when $\|\Omega\|_{L^2}$ is arbitrarily small. Smooth truncations further show that the corresponding uniform interior $W^{2,1}$ estimate fails even for smooth solutions and coefficients. Thus, the $L\log L$ assumption ensuring local $W^{2,1}$ regularity cannot in general be replaced by componentwise local Hardy data.
\end{abstract}

\maketitle

\maketitle
	
	{\small
		\smallskip
		\newline
	}
	\bigskip
	
	\arraycolsep=1pt

\section{Introduction and main results}\label{sec:introduction}

In this paper, we consider the regularity of weak solutions
$u\in W^{1,2}(B_1,\R^n)$ to the following second order elliptic system of
Rivi\`ere type:
\begin{equation}\label{eq:system}
-\Delta u=\Omega\cdot\nabla u+f\qquad\text{in }B_1,
\end{equation}
where $B_\rho=\{x\in\R^2:|x|<\rho\}\subset \R^2$ and
$\Omega=(\Omega_1,\Omega_2)\in L^2(B_1,so(n)\otimes\R^2)$.
Here $so(n)$ denotes the space of $n\times n$ antisymmetric matrices. In local coordinates, \eqref{eq:system} reads as
\[
-\Delta u^i=\Omega^i_j\cdot \nabla u^j+f^i\qquad \text{for all }i=1,\cdots,n. 
\]

The homogeneous version of \eqref{eq:system}, that is when $f\equiv0$,
was introduced by Rivi\`ere in his celebrated work \cite{Riviere} on conservation law for critical points of conformally invariant elliptic Lagrangian. Since $\Omega\cdot\nabla u$ belongs merely to $L^1$, the standard elliptic
regularity theory does not directly give continuity of weak solutions.
Rivi\`ere's key observation was that the antisymmetry of $\Omega$ allows
one to find a bounded Sobolev gauge and rewrite the system in conservation
law form. This leads to continuity of weak solutions under the natural
$L^2$ assumption on the potential. For a comprehensive introduction on the regularity theory of Rivi\`ere's equation and its  applications, see \cite{Riviere-2012, Lamm-2007, LammSharp2016}. 

For the inhomogeneous system, Sharp--Topping \cite{SharpTopping}
established local $W^{2,p}$ regularity when $f\in L^p$ and $1<p<2$.
They also treated the endpoint case $f\in L\log L$ and obtained local
$W^{2,1}$ regularity; see \cite[Theorems 1.1 and 1.6]{SharpTopping}. This result was extended by Du--Kang--Wang \cite{Du-Kang-Wang-2024-PAMS} to the case $f$ belongs to certain Morrey spaces $M^{p,s}$. 
Related $L^p$, $L\log L$ and $M^{p,s}$ regularity results for higher order system of Rivi\`ere type were developed in \cite{GXZ,GXZ2022,Guo-Wang-Xiang-2023-CVPDE,Guo-Qi-2025,Xiang-Zheng-2024}. 

It is therefore natural to ask what happens in the critical case, namely, \emph{whether $f$ belongs to the local Hardy space $h^1_{\mathrm{loc}}$ ensures local $W^{2,1}$ regularity for \eqref{eq:system}}. A related question
was raised by Lamm--Sharp \cite[Remark~2.2(5)]{LammSharp2016} in their discussion
of estimates with an additional control on the Hopf differential.
Here we consider \eqref{eq:system} without such a condition; see Remark \ref{thm:hopf}.

Note that if the drift term $f\in h^1_{\mathrm{loc}}$, then the strong $W^{1,2}$-compactness for \eqref{eq:system} fails by \cite[Section 8]{SharpTopping}, while as contrast the strong $W^{1,2}$-compactness do hold for $f\in L\log L$ \cite{Li-Zhu-2012}. This indicates that there is an essential difference between the case $f\in L\log L$ and $f\in h^1_{\mathrm{loc}}$. On the other hand, in the example of Sharp--Topping \cite[Section 8]{SharpTopping}, $\Omega=0$, and so each individual solution still
belongs to $W_{\mathrm{loc}}^{2,1}$ by the Calder\'on--Zygmund theory in local Hardy spaces; see \cite[Appendix A.2]{SharpTopping}. Thus, it does not settle the above critical regularity question for a single
solution of \eqref{eq:system}.

Throughout the paper, vector-valued
functions in Hardy space are always understood in the componentwise sense, and the norm is the sum
of the corresponding component norms. In particular, $f\in h^1$ does not require
$|f|\in h^1$. 

Our main result of this paper shows that the assumption $f\in h^1_{\mathrm{loc}} $ does not
ensure local $W^{2,1}$ regularity, even when $\|\Omega\|_{L^2}$
is arbitrarily small.

\begin{theorem}\label{thm:main}
For every $n\ge2$ and every $\varepsilon>0$, there exist
\begin{align*}
&\Omega\in L^2(B_1,so(n)\otimes\R^2),\qquad
  \|\Omega\|_{L^2(B_1)}<\varepsilon,\\
&u\in W^{1,2}_0(B_1,\R^n)\cap C(\overline{B_1},\R^n),\qquad
  f\in h^1(B_1,\R^n),
\end{align*}
such that \eqref{eq:system} holds in the sense of distributions, but
\begin{equation}\label{eq:nonw21}
u\notin W^{2,1}(B_\rho,\R^n)\qquad\text{for every }0<\rho<1.
\end{equation}
Moreover, the zero extension $f^0$ belongs to $h^1(\R^2,\R^n)$, and
there exists $\widetilde f\in\mathcal H^1(\R^2,\R^n)$ which agrees with
$f$ in $B_1$. The quantities
$\|\Omega\|_{L^2(B_1)}$, $\|u\|_{W^{1,2}(B_1)}$,
$\|u\|_{L^\infty(B_1)}$ and $\|f^0\|_{h^1(\R^2)}$
can be made arbitrarily small simultaneously.
\end{theorem}


The next result shows that the corresponding uniform interior
estimate fails even for smooth solutions and coefficients.

\begin{theorem}\label{thm:smooth}
For every $\varepsilon>0$, there exist sequences
\begin{align*}
u_N&\in C^\infty(\overline{B_1},\R^2),\qquad
f_N\in C^\infty(\overline{B_1},\R^2),\\
\Omega_N&\in C^\infty(\overline{B_1},so(2)\otimes\R^2)
\end{align*}
such that
\[
-\Delta u_N=\Omega_N\cdot\nabla u_N+f_N\quad\text{in }B_1,
\qquad u_N|_{\partial B_1}=0,
\qquad \|\Omega_N\|_{L^2(B_1)}<\varepsilon,
\]
and
\[
\sup_N\left(\|u_N\|_{W^{1,2}(B_1)}+\|u_N\|_{L^\infty(B_1)}
+\|f_N^0\|_{h^1(\R^2)}\right)<\infty,
\]
where $f_N^0$ is the zero extension of $f_N$, whereas
\[
\|D^2u_N\|_{L^1(B_{1/2})}\to\infty.
\]
\end{theorem}

As an immediate consequence, we obtain the following failure of a
uniform interior estimate.

\begin{corollary}\label{cor:noestimate}
There are no positive constants $\eta$ and $C$ such that every smooth
solution of \eqref{eq:system} with $\|\Omega\|_{L^2(B_1)}\le\eta$
satisfies
\[
\|u\|_{W^{2,1}(B_{1/2})}
\le C\left(\|u\|_{L^1(B_1)}+\|f\|_{h^1(B_1)}\right).
\]
The conclusion remains valid if $\|u\|_{L^1(B_1)}$ is replaced by
$\|u\|_{W^{1,2}(B_1)}$, or if $\|f\|_{h^1(B_1)}$ is replaced by
$\|f^0\|_{h^1(\R^2)}$.
\end{corollary}

\begin{proof}
Apply Theorem~\ref{thm:smooth} with $\varepsilon<\eta$. The quantities
on the right-hand side remain bounded, while
$\|u_N\|_{W^{2,1}(B_{1/2})}\ge\|D^2u_N\|_{L^1(B_{1/2})}\to\infty$.
For the intrinsic Hardy norm, we use
$\|f_N\|_{h^1(B_1)}\le\|f_N^0\|_{h^1(\R^2)}$.
\end{proof}

The construction uses a radial map whose derivative with respect
to logarithmic radius has a slowly rotating direction.
The potential is chosen so that the term generated by this
rotation is absorbed into $\Omega\cdot\nabla u$.
An annular atomic decomposition shows that the remaining forcing
belongs to local Hardy space, whereas a tangential component of
the Hessian is not integrable.
An explicit orthogonal gauge transforms the forcing into a fixed
vector multiplied by a nonnegative scalar function outside local
Hardy space.
Thus, multiplication by a bounded $W^{1,2}$ gauge, which preserves
$L\log L$, need not preserve componentwise local Hardy space
membership.

This paper is organized as follows. Section~\ref{sec:hardy} contains
some preliminaries and auxiliary estimates. We prove
Theorem~\ref{thm:main} in Section~\ref{sec:construction} and discuss
the associated gauge in Section~\ref{sec:gauge}. The proof of
Theorem~\ref{thm:smooth} is given in Section~\ref{sec:smooth}.
Throughout the paper, $C$ denotes a positive constant which may change
from line to line. Its dependence on the dimension or a fixed cutoff is
allowed, but it is independent of the parameters $T$ and $N$ unless
otherwise stated.

\section{Preliminaries and auxiliary results}\label{sec:hardy}

In this section, we collect the Hardy space estimates used in the
construction. We first recall the maximal function definition and an
elementary estimate for atoms. We then give an annular decomposition
which will be applied to the forcing term. For a detailed introduction on the theory of Hardy spaces, one can see \cite{Stein-1993,Stein-Weiss-1971}. Here we mainly follow the presentation in \cite{Lamm-2007, Goldberg, SharpTopping}. 

\subsection{Hardy spaces and atoms}
Fix a nonnegative function $\varphi\in C_c^\infty(B_1)$ with
$\int_{\R^2}\varphi=1$ which is positive on $\overline{B_{1/2}}$,
and write $\varphi_s(x)=s^{-2}\varphi(x/s)$. For a locally integrable
function $F$ on $\R^2$, set
\[
\mathcal M F(x)=\sup_{s>0}|\varphi_s*F(x)|,
\qquad
\mathcal M_{\mathrm{loc}}F(x)=\sup_{0<s<1}|\varphi_s*F(x)|.
\]
The corresponding $L^1$ norms give equivalent norms on
$\HH(\R^2)$ and Goldberg's local Hardy space $\hh(\R^2)$,
respectively; see \cite{Goldberg}. On $B_1$, we use the intrinsic
maximal function with $0<s<1-|x|$, as in
\cite[Appendix A.2]{SharpTopping}. In the construction below, we shall
prove the stronger statement that the zero extension of the forcing
belongs to $\hh(\R^2)$.

Recall that a bounded function $a$ supported in a ball $B=B(x_0,r)$
is an $\HH$ atom if
\[
\int_{\R^2}a=0,
\qquad \|a\|_{L^\infty}\le |B|^{-1}.
\]
Such an atom satisfies $\|a\|_{\HH}\le C$, with a constant
independent of $B$. For completeness, we recall the short proof.
On $2B$, we have
$\mathcal Ma\le\|\varphi\|_{L^1}\|a\|_{L^\infty}$, and hence
$\int_{2B}\mathcal Ma\le C$. If $d=|x-x_0|>2r$, a nonzero
convolution requires $s\ge d/2$. By the cancellation of $a$ and the
mean value theorem,
\begin{align*}
|\varphi_s*a(x)|
&=\left|\int_B\bigl(\varphi_s(x-y)-\varphi_s(x-x_0)\bigr)a(y)\,dy\right|\\
&\le Crs^{-3}\|a\|_{L^1}
\le Crd^{-3}.
\end{align*}
Consequently,
\[
\int_{\R^2\setminus2B}\mathcal Ma(x)\,dx
\le Cr\int_{2r}^\infty s^{-2}\,ds\le C,
\]
which proves the assertion. It follows by subadditivity of the maximal
function and Fatou's lemma that an absolutely summable linear combination
of atoms belongs to $\HH$. The same argument applies componentwise
to vector-valued functions, with a constant depending on the target
dimension.

\subsection{An annular decomposition}
The following lemma gives a convenient sufficient condition for Hardy space membership of the radial forcing used below.

\begin{lemma}\label{lem:annular}
Let $T\ge2$ and let $F\colon [T,\infty)\to\R^d$ be measurable.
Suppose that, for some $K>0$,
\begin{equation}\label{eq:tail-assumptions}
|F(t)|\le Kt^{-2},\qquad
\left|\int_t^\infty F(s)\,ds\right|\le Kt^{-4/3}
\quad\text{for }t\ge T.
\end{equation}
Define
\[
g(x)=|x|^{-2}F\bigl(T+\log(1/|x|)\bigr),\qquad 0<|x|<1.
\]
Then the zero extension $g^0$ belongs to $\hh(\R^2,\R^d)$.
Moreover, there exists $\widetilde g\in\HH(\R^2,\R^d)$ which
agrees with $g$ in $B_1$, and
\begin{equation}\label{eq:hardybound}
\|g^0\|_{\hh(\R^2)}+\|\widetilde g\|_{\HH(\R^2)}
\le C_dKT^{-1/3}.
\end{equation}
\end{lemma}

\begin{proof}
By the definition of \(g\) and polar coordinates,
\[
\begin{aligned}
\|g\|_{L^1(B_1)}
&=2\pi\int_0^1 r^{-2}
\left|F\!\left(T+\log\frac1r\right)\right|r\,dr\\
&=2\pi\int_T^\infty |F(t)|\,dt
\leq\frac{2\pi K}{T}.
\end{aligned}
\]
For $k=0,1,\ldots$, put $r_k=e^{-k}$, $B_k=B_{r_k}$ and
$D_k=B_k\setminus B_{k+1}$. Set
\[
\beta_k=|B_k|^{-1}\mathbf1_{B_k},\qquad
g_k=g\mathbf1_{D_k},\qquad m_k=\int_{\R^2}g_k.
\]
Using $t=T+\log(1/r)$, we obtain
\[
m_k=2\pi\int_{T+k}^{T+k+1}F(t)\,dt,
\qquad
S_k:=\sum_{j=k}^\infty m_j=2\pi\int_{T+k}^\infty F(t)\,dt.
\]
It follows from \eqref{eq:tail-assumptions} that
\begin{equation}\label{eq:annular-bounds}
|m_k|\le CK(T+k)^{-2},\qquad
|S_k|\le CK(T+k)^{-4/3}.
\end{equation}

We first remove the mean on each annulus. The function
$q_k=g_k-m_k\beta_k$ is supported in $B_k$ and has integral zero.
For $x\in D_k$, write $r=|x|$ and $t=T+\log(1/r)$. Then
$t\ge T+k$ and
$r^{-2}\le e^{2(k+1)}=\pi e^2/|B_k|$. Hence, by
\eqref{eq:tail-assumptions}, \eqref{eq:annular-bounds},
and $\|\beta_k\|_{L^\infty}=|B_k|^{-1}$,
\begin{align*}
\|q_k\|_{L^\infty}
     \le \frac{\pi e^2K(T+k)^{-2}+|m_k|}{|B_k|}
\le \frac{CK(T+k)^{-2}}{|B_k|}.
\end{align*}
Thus, $\|q_k\|_{\HH}\le C_dK(T+k)^{-2}$.

Next, we estimate the remaining means. Since $m_k=S_k-S_{k+1}$,
we can rewrite the sum as
\[
\sum_{k=0}^{M}m_k\beta_k
=S_0\beta_0+\sum_{k=1}^{M}S_k(\beta_k-\beta_{k-1})
-S_{M+1}\beta_M.
\]
The last term tends to zero in $L^1$ by \eqref{eq:annular-bounds}.
Hence
\begin{equation*}\label{eq:telescoping}
\sum_{k=0}^\infty m_k\beta_k
=S_0\beta_0+\sum_{k=1}^\infty S_k(\beta_k-\beta_{k-1})
\end{equation*}
in $L^1$. Each difference $\beta_k-\beta_{k-1}$ has integral zero,
is supported in $B_{k-1}$, and satisfies
\[
\|\beta_k-\beta_{k-1}\|_{L^\infty}
\le\frac{e^2+1}{|B_{k-1}|}.
\]
Therefore these differences are fixed multiples of atoms.

To obtain a global Hardy extension, choose a fixed smooth function
$\psi$ supported in $B_2\setminus\overline{B_1}$ with integral one,
and define $\widetilde g=g^0-S_0\psi$. Combining the preceding
decompositions gives
\begin{equation*}\label{eq:full-atoms}
\widetilde g
=\sum_{k=0}^\infty q_k
+\sum_{k=1}^\infty S_k(\beta_k-\beta_{k-1})
+S_0(\beta_0-\psi).
\end{equation*}
The function $\beta_0-\psi$ is bounded, supported in $B_2$, and has
integral zero. Applying the atomic estimate and
\eqref{eq:annular-bounds}, we conclude that
\begin{align*}
\|\widetilde g\|_{\HH}
&\le C_dK\left(\sum_{k=0}^\infty(T+k)^{-2}
+\sum_{k=1}^\infty(T+k)^{-4/3}+T^{-4/3}\right)\\
&\le C_dK(T^{-1}+T^{-1/3})
\le C_dKT^{-1/3}.
\end{align*}
Finally, $\HH\hookrightarrow\hh$ since
$\mathcal M_{\mathrm{loc}}\le\mathcal M$.
The fixed bounded function $\psi$ also belongs to $\hh$: its local
maximal function is bounded and supported in $B_3$.
Using $g^0=\widetilde g+S_0\psi$, we obtain the same bound for
$\|g^0\|_{\hh}$. This proves the lemma.
\end{proof}

\subsection{Oscillatory integral estimates}
For $t\ge2$, define
\begin{equation}\label{eq:phase}
\theta(t)=3t^{1/3},\qquad
E(t)=\begin{pmatrix}\cos\theta(t)\\ \sin\theta(t)\end{pmatrix},\qquad
J=\begin{pmatrix}0&-1\\1&0\end{pmatrix}.
\end{equation}
Then $\theta'(t)=t^{-2/3}$ and $E'(t)=t^{-2/3}JE(t)$.
We shall use the following elementary estimates.

\begin{lemma}\label{lem:tails}
For every $t\ge2$, the improper integrals below converge and satisfy
\begin{equation}\label{eq:two-tails}
\left|\int_t^\infty\frac{E(s)}s\,ds\right| \le 2t^{-1/3},
\qquad
\left|\int_t^\infty\frac{E(s)}{s^2}\,ds\right| \le 2t^{-4/3}.
\end{equation}
Moreover, if $b\in C^1([T,\infty))$, $T\ge2$, satisfies
\begin{equation}\label{eq:envelope}
|b(t)|\le Kt^{-2},\qquad |b'(t)|\le Kt^{-3},
\end{equation}
then
\begin{equation}\label{eq:uniform-tail}
\left|\int_t^\infty b(s)E(s)\,ds\right|\le CKt^{-4/3}
\qquad\text{for }t\ge T.
\end{equation}
\end{lemma}

\begin{proof}
Identifying $\R^2$ with $\mathbb C$, we write $E(s)=e^{i\theta(s)}$.
Integration by parts gives
\begin{align*}
\int_t^R b(s)e^{i\theta(s)}\,ds
=\left[\frac{b(s)}{i\theta'(s)}e^{i\theta(s)}\right]_t^R
-\int_t^R\left(\frac{b(s)}{i\theta'(s)}\right)' e^{i\theta(s)}\,ds.
\end{align*}
For $b(s)=s^{-1}$, then $b(s)/\theta'(s)$ equals
$s^{-1/3}$. Its boundary value at infinity vanishes, and its derivative
is integrable. Letting $R\to\infty$, we obtain
\[
\left|\int_t^\infty\frac{E(s)}s\,ds\right|
\le t^{-1/3}+\frac13\int_t^\infty s^{-4/3}\,ds
=2t^{-1/3}.
\]
The same argument with $b(s)=s^{-2}$ gives
\[
\left|\int_t^\infty\frac{E(s)}{s^2}\,ds\right|
\le t^{-4/3}+\frac43\int_t^\infty s^{-7/3}\,ds
=2t^{-4/3}.
\]
In particular, the first integral converges despite the lack of absolute
integrability of $E(s)/s$.

Finally, under \eqref{eq:envelope},
\[
|b(t)t^{2/3}|\le Kt^{-4/3},\qquad
|(b(t)t^{2/3})'|\le CKt^{-7/3}.
\]
The same integration by parts, followed by integration of the last
bound, proves \eqref{eq:uniform-tail}.
\end{proof}

\section{Proof of Theorem~\ref{thm:main}}\label{sec:construction}

This section is devoted to the proof of Theorem~\ref{thm:main}.
We first construct the example with two target components and verify
its regularity and the equation. We then apply Lemma~\ref{lem:annular}
to the forcing and estimate the tangential component of the Hessian.

\begin{proof}[Proof of Theorem~\ref{thm:main}]
It suffices to consider $n=2$. For $n>2$, we add zero components to
$u$ and $f$, and extend each $\Omega_j$ by a zero block. This preserves
antisymmetry, the equation and all the estimates below.
We divide the proof into four steps.

\medskip
\noindent\textbf{Step 1. Construction and elementary estimates.}\par\nobreak

Fix $T\ge2$. With $E$ and $J$ as in \eqref{eq:phase}, set
\begin{equation}\label{eq:tU}
t=T+\log(1/r),\qquad
U(t)=-\int_t^\infty\frac{E(s)}s\,ds,\qquad
u_T(x)=U(t)-U(T),
\end{equation}
where $r=|x|\in(0,1]$, and define $u_T(0)=-U(T)$.
For $x\ne0$, let
\begin{equation}\label{eq:Omegaf}
\Omega_{T,j}(x)=\frac{x_j}{r^2t^{2/3}}J\quad(j=1,2),
\qquad
f_T(x)=\frac{E(t)}{r^2t^2}.
\end{equation}
The values of $\Omega_T$ and $f_T$ at the origin may be chosen
arbitrarily. Since $J$ is antisymmetric, so is each $\Omega_{T,j}$.

By \eqref{eq:two-tails}, $U(t)\to0$ as $t\to\infty$.
It follows that $u_T\in C(\overline{B_1},\R^2)$ and
$u_T=0$ on $\partial B_1$. Moreover,
\begin{equation}\label{eq:sup-bound}
\|u_T\|_{L^\infty(B_1)}\le4T^{-1/3}.
\end{equation}
Write $V(t)=U'(t)=t^{-1}E(t)$. Since $\partial_jt=-x_j/r^2$,
we have, away from the origin,
\[
\partial_j u_T=-\frac{x_j}{r^2t}E(t),\qquad
|\nabla u_T|=\frac1{rt}.
\]
Using polar coordinates and $dr/r=-dt$, we obtain
\begin{align}
\|\nabla u_T\|_{L^2(B_1)}^2
&=2\pi\int_T^\infty t^{-2}\,dt
=\frac{2\pi}{T},\label{eq:energy}\\
\|\Omega_T\|_{L^2(B_1)}^2
&=4\pi\int_T^\infty t^{-4/3}\,dt
=12\pi T^{-1/3},\label{eq:potential-energy}\\
\|f_T\|_{L^1(B_1)}
&=2\pi\int_T^\infty t^{-2}\,dt
=\frac{2\pi}{T}.\label{eq:fL1}
\end{align}
Here the factor $4\pi$ in \eqref{eq:potential-energy} comes from
$|J|^2=2$ and $\sum_{j=1}^2x_j^2=r^2$.

Let $\varphi\in C_c^\infty(B_1)$ and $j=1,2$. Integration by parts on $B_1\setminus\overline{B_\delta}$ gives
\begin{align*}
\int_{B_1\setminus\overline{B_\delta}} u_T\,\partial_j\varphi\,dx
=\int_{B_1\setminus\overline{B_\delta}} \frac{x_j}{r^2t}E(t)\varphi\,dx
-\int_{\partial B_\delta} u_T\varphi\,\frac{x_j}{\delta}\,d\sigma.
\end{align*}
Then, the boundary term satisfies
\begin{align*}
\left| \int_{\partial B_\delta} u_T\varphi\,\frac{x_j}{\delta}\,d\sigma \right|
\le 2\pi\delta\, \|u_T\|_{L^\infty(B_1)} \|\varphi\|_{L^\infty(B_1)}
\to 0
\end{align*}
as $\delta\to0$. Letting $\delta\to0$ and using \eqref{eq:sup-bound}–\eqref{eq:energy}, we obtain
\[
\int_{B_1}u_T\,\partial_j\varphi\,dx
=
\int_{B_1}\frac{x_j}{r^2t}E(t)\varphi\,dx.
\]
Thus $u_T\in W^{1,2}(B_1,\mathbb R^2)$. Its zero boundary trace yields $u_T\in W^{1,2}_0(B_1,\mathbb R^2).$

\medskip
\noindent\textbf{Step 2. Verification of the weak equation.}\par\nobreak

For the radial function in \eqref{eq:tU}, direct differentiation gives
\[
\partial_r u_T=-\frac{V(t)}r,\qquad
\partial_r^2u_T=\frac{V'(t)+V(t)}{r^2},\qquad
\Delta u_T=\partial_r^2u_T+\frac1r\partial_r u_T
=\frac{V'(t)}{r^2}.
\]
By the definition of $E$,
\begin{equation}\label{eq:Vprime}
V'(t)=-t^{-2}E(t)+t^{-5/3}JE(t).
\end{equation}
On the other hand,
\begin{equation}\label{eq:drift}
\Omega_T\cdot\nabla u_T
=-\sum_{j=1}^2\frac{x_j^2}{r^4t^{5/3}}JE(t)
=-\frac{JE(t)}{r^2t^{5/3}}.
\end{equation}
Combining \eqref{eq:Omegaf}, \eqref{eq:Vprime} and \eqref{eq:drift},
we obtain
\[
-\Delta u_T=\Omega_T\cdot\nabla u_T+f_T
\qquad\text{in }B_1\setminus\{0\}.
\]
Both terms on the right-hand side are integrable. Indeed,
\[
\|\Omega_T\cdot\nabla u_T\|_{L^1(B_1)}
=2\pi\int_T^\infty t^{-5/3}\,dt
=3\pi T^{-2/3},
\]
and the estimate for $f_T$ was given in \eqref{eq:fL1}.

It remains to verify the equation at the origin. For
$\zeta\in C_c^\infty(B_1,\R^2)$, integration by parts on
$B_1\setminus\overline{B_\delta}$ gives
\begin{align*}
\int_{B_1\setminus B_\delta}\nabla u_T\cdot\nabla\zeta
&=\int_{B_1\setminus B_\delta}
(\Omega_T\cdot\nabla u_T+f_T)\cdot\zeta
-\int_{\partial B_\delta}\partial_r u_T\cdot\zeta\,d\sigma.
\end{align*}
Furthermore,
\begin{equation}\label{eq:flux}
\left|\int_{\partial B_\delta}\partial_r u_T\cdot\zeta\,d\sigma\right|
\le\frac{2\pi\|\zeta\|_{L^\infty}}
{T+\log(1/\delta)}\to0.
\end{equation}
Passing to the limit, using $\nabla u_T\in L^2$ and the integrability
of the right-hand side, we conclude that
\[
\int_{B_1}\nabla u_T\cdot\nabla\zeta
=\int_{B_1}(\Omega_T\cdot\nabla u_T+f_T)\cdot\zeta.
\]
Thus \eqref{eq:system} holds on the whole disk in the sense of
distributions.

\medskip
\noindent\textbf{Step 3. Hardy space estimate for the forcing.}\par\nobreak

We apply Lemma~\ref{lem:annular} with $F(t)=t^{-2}E(t)$.
The pointwise bound in \eqref{eq:tail-assumptions} follows from
$|E(t)|=1$, and the second estimate in \eqref{eq:two-tails} gives the required integral bound. Hence
\begin{equation}\label{eq:forcinghardy}
f_T^0\in\hh(\R^2,\R^2),\qquad
\|f_T^0\|_{\hh(\R^2)}\le CT^{-1/3}.
\end{equation}
The same lemma provides an extension
$\widetilde f_T\in\HH(\R^2,\R^2)$ satisfying
$\widetilde f_T=f_T$ in $B_1$ and
$\|\widetilde f_T\|_{\HH}\le CT^{-1/3}$.
In particular, $f_T\in h^1(B_1,\R^2)$.

\medskip
\noindent\textbf{Step 4. Failure of $W^{2,1}$ regularity.}\par\nobreak

Since $u_T$ is radial, for $x\ne0$ we have
\[
|D^2u_T|^2
=
|\partial_r^2u_T|^2
+\frac{|\partial_ru_T|^2}{r^2}
\ge\frac{1}{r^4t^2}.
\]
Thus $|D^2u_T|\ge 1/(r^2t)$, and for every $0<\rho<1$,
\begin{equation}\label{eq:divergence}
\begin{aligned}
\int_{B_\rho\setminus\{0\}}|D^2u_T|\,dx
&\ge
2\pi\int_0^\rho
\frac{dr}{r(T+\log(1/r))}\\
&=
2\pi\int_{T+\log(1/\rho)}^\infty\frac{dt}{t}
=\infty.
\end{aligned}
\end{equation}
If $u_T\in W^{2,1}(B_\rho)$, its weak Hessian would agree
with the classical Hessian on $B_\rho\setminus\{0\}$,
contradicting \eqref{eq:divergence}.
This proves \eqref{eq:nonw21}.

Finally, \eqref{eq:sup-bound}, \eqref{eq:energy},
\eqref{eq:potential-energy} and \eqref{eq:forcinghardy} imply
\[
\|\Omega_T\|_{L^2(B_1)}
+\|u_T\|_{W^{1,2}(B_1)}
+\|u_T\|_{L^\infty(B_1)}
+\|f_T^0\|_{h^1(\mathbb R^2)}
\to0
\qquad\text{as }T\to\infty.
\]
Choosing $T$ sufficiently large gives all the smallness
assertions and completes the proof.
\end{proof}

\begin{remark}\label{rem:orlicz}
The forcing $f_T$ does not belong to $L\log L$ near the origin.
Indeed, $|f_T|=r^{-2}t^{-2}$ and
$\log|f_T|=2(t-T)-2\log t$. Thus
$\log(2+|f_T|)\ge ct$ for all sufficiently large $t$.
For every $0<\rho<1$, we may therefore choose $t_0$ sufficiently
large, depending on $T$ and $\rho$, such that
\[
\int_{B_\rho}|f_T|\log(2+|f_T|)\,dx
\ge c\int_{t_0}^\infty\frac{dt}{t}=\infty.
\]
Hence the example is consistent with the endpoint $L\log L$ result
of Sharp--Topping \cite[Theorem 1.6]{SharpTopping}.
\end{remark}

\section{The associated gauge and Hardy space cancellation}\label{sec:gauge}

In this section, we compute an explicit gauge for the example in
Section~\ref{sec:construction}. The calculation shows that multiplication
by a bounded $W^{1,2}$ gauge can remove the cancellation responsible
for the Hardy space membership of the forcing.

With the same notation as above, define
\begin{equation*}\label{eq:gauge}
A_T(x)=\exp\bigl(-\theta(t)J\bigr)\in\mathrm{SO}(2),
\qquad x\ne0.
\end{equation*}
Then $A_TE(t)=e_1$, where $e_1=(1,0)^\mathsf T$.
Since $J$ commutes with its exponential, we have
\begin{equation}\label{eq:puregauge}
\partial_jA_T
=-\theta'(t)JA_T\partial_jt
=\frac{x_j}{r^2t^{2/3}}A_TJ
=A_T\Omega_{T,j}.
\end{equation}
Orthogonality gives $|\partial_jA_T|=|\Omega_{T,j}|$ and hence
\[
A_T\in L^\infty(B_1)\cap W^{1,2}(B_1),\qquad
\|\nabla A_T\|_{L^2(B_1)}=\|\Omega_T\|_{L^2(B_1)}.
\]
To verify the weak derivatives at the origin, we use the same integration
by parts argument as in Step~1 of the preceding proof. The boundary
term is bounded by $C\delta\|A_T\|_{L^\infty}$ times the supremum
norm of the test function, and therefore tends to zero.
Thus \eqref{eq:puregauge} is Rivi\`ere's gauge relation
$\nabla A_T-A_T\Omega_T=\nabla^\perp B$ with $B=0$.

For this gauge, the transformed forcing is
\begin{equation*}\label{eq:positiveg}
A_Tf_T=g_Te_1,\qquad
g_T(x)=\frac1{r^2t^2}\ge0.
\end{equation*}
For \(x\ne0\), the identity \(A_TE(t)=e_1\) gives
\[
A_T\partial_j u_T
=-\frac{x_j}{r^2t}\,e_1,
\qquad j=1,2.
\]
Consequently,
\[
\operatorname{div}(A_T\nabla u_T)
=\sum_{j=1}^2\partial_j(A_T\partial_j u_T)
=-\frac{1}{r^2t^2}e_1
=-g_Te_1.
\]
The latter identity holds on the whole disk in the sense of distributions.
Indeed, the normal flux across $\partial B_\delta$ has magnitude
$2\pi/(T+\log(1/\delta))\to0$, so no additional term appears at
the origin.

\begin{proposition}\label{prop:multiplier}
For every $0<\rho<1$, one has $g_T\notin h^1(B_\rho)$.
Consequently, $A_Tf_T$ does not belong to local Hardy space near the
origin, although $f_T$ does, $A_T$ is bounded and orthogonal, and
$\|\nabla A_T\|_{L^2(B_1)}$ can be arbitrarily small.
\end{proposition}

\begin{proof}
By polar integration,
\[
\int_{B_r}g_T\,dx
=2\pi\int_{T+\log(1/r)}^\infty t^{-2}\,dt
=\frac{2\pi}{T+\log(1/r)}.
\]
Fix $x$ with $|x|=r$. For $y\in B_r$, we have $|x-y|\le2r$.
Since $\varphi$ is bounded below by a positive constant on
$\overline{B_{1/2}}$, it follows that
\[
(\varphi_{4r}*g_T)(x)
\ge\frac{c}{r^2}\int_{B_r}g_T(y)\,dy
\ge\frac{c}{r^2(T+\log(1/r))}.
\]
Here, for $r<\rho/5$, the convolution only uses values in $B_\rho$
and $4r<\rho-|x|$, so this scale is admissible in the intrinsic
maximal function on $B_\rho$. Integrating over $B_{\rho/5}$, we obtain
\[
\int_{B_{\rho/5}}\sup_{0<s<\rho-|x|}
|\varphi_s*g_T(x)|\,dx
\ge c\int_0^{\rho/5}\frac{dr}{r(T+\log(1/r))}
=\infty.
\]
This proves the proposition.
\end{proof}

It follows that the implication
\[
f\in h^1,\quad A\in W^{1,2}\cap L^\infty
\quad\to\quad Af\in h^1
\]
fails even for the orthogonal gauges constructed above. In the
$L\log L$ argument, multiplication by a bounded gauge preserves the
class of the forcing. Proposition~\ref{prop:multiplier} shows that this
step is unavailable for general componentwise Hardy data, and
Theorem~\ref{thm:main} gives the corresponding failure of regularity
for a weak solution.

\begin{remark}\label{rem:scope}
We record two features which specify the scope of the example.
First, $|f_T|=g_T$, so Proposition~\ref{prop:multiplier} also shows
that $|f_T|$ does not belong to local Hardy space near the origin.
Thus the example does not apply to a hypothesis imposed on the Hardy
maximal function of $|f|$, nor to one which additionally requires a
suitable gauge-transformed forcing to belong to Hardy space.

Second, the example does not satisfy a fixed growth bound of the form
$|\Omega|\le C|\nabla u|$. In fact,
\[
\frac{|\Omega_T(x)|}{|\nabla u_T(x)|}
=\sqrt2\,t^{1/3}\to\infty
\qquad\text{as }x\to0.
\]
Accordingly, no conclusion is drawn here for classes of approximate
harmonic maps with additional target or growth assumptions. The
construction also concerns only the second order system
\eqref{eq:system}, and gives no conclusion for fourth order
Lamm--Rivi\`ere systems.
\end{remark}

\begin{remark}\label{thm:hopf}
The estimate of Lamm and Sharp
\cite[Theorem~2.1]{LammSharp2016} also assumes that
$|\Phi|^{1/2}\in L_{\mathrm{loc}}^{2,1}$, where
\[
\Phi
=
|\partial_1u|^2-|\partial_2u|^2
-2i\langle\partial_1u,\partial_2u\rangle
\]
is the Hopf differential.
For the radial map $u_T$ constructed above, we have
\[
|\Phi_T|^{1/2}
=
|\partial_r u_T|
=
\frac{1}{r\bigl(T+\log(1/r)\bigr)}.
\]
The radial characterization of $L^{2,1}$, together with
\[
\int_0^\rho
\frac{dr}{r\bigl(T+\log(1/r)\bigr)}
=
\infty,
\qquad 0<\rho<1,
\]
shows that $|\Phi_T|^{1/2}\notin L^{2,1}(B_\rho)$.
Thus our construction here does not automatically solve the open question of Lamm--Sharp. 
\end{remark}

\section{Proof of Theorem~\ref{thm:smooth}}\label{sec:smooth}

In this section, we truncate the construction from
Section~\ref{sec:construction} to obtain smooth solutions. The cutoff
is chosen on the logarithmic scale so that the estimates in
Lemma~\ref{lem:tails} remain uniform. This gives a uniform Hardy space
bound for the forcing, while the Hessian norm diverges as the cutoff
is removed.

\begin{proof}[Proof of Theorem~\ref{thm:smooth}]
Fix $T\ge2$ sufficiently large that
$12\pi T^{-1/3}<\varepsilon^2$. Choose
$\chi\in C^\infty([0,\infty))$ satisfying $0\le\chi\le1$,
$\chi=1$ on $[0,1]$ and $\chi=0$ on $[2,\infty)$.
We divide the proof into three steps.

\medskip
\noindent\textbf{Step 1. Smooth truncation and the equation.}\par\nobreak

For integers $N>2T$, define
\begin{align*}
V_N(t)&=\frac{\chi(t/N)}t E(t),\qquad
U_N(t)=-\int_t^\infty V_N(s)\,ds,\\
u_N(x)&=U_N(t)-U_N(T),\qquad
\Omega_{N,j}(x)=\chi\bigl(t/(2N)\bigr)\Omega_{T,j}(x),
\end{align*}
where $t=T+\log(1/r)$ and $r=|x|>0$. Set
\begin{equation}\label{eq:bn}
b_N(t)=\frac{\chi(t/N)}{t^2}-\frac{\chi'(t/N)}{Nt},
\qquad
f_N(x)=\frac{b_N(t)}{r^2}E(t).
\end{equation}
For $t\ge2N$, the function $u_N$ is constant and $f_N=0$.
Moreover, $\Omega_N=0$ for $t\ge4N$. Defining these functions by
their corresponding constant values at the origin, we see that all
three are smooth on $\overline{B_1}$. Smoothness across the cutoff
boundaries follows from the choice of $\chi$. Also,
$u_N=0$ on $\partial B_1$.

Differentiating $V_N$, we have
\[
-V_N'(t)=b_N(t)E(t)-\frac{\chi(t/N)}{t^{5/3}}JE(t).
\]
Whenever $\chi(t/N)$ is nonzero, we have $t<2N$ and hence
$\chi(t/(2N))=1$. Therefore
\[
\Omega_N\cdot\nabla u_N
=-\frac{\chi(t/N)}{r^2t^{5/3}}JE(t).
\]
Since $\Delta u_N=r^{-2}V_N'(t)$, these identities prove
$-\Delta u_N=\Omega_N\cdot\nabla u_N+f_N$ away from the origin.
The same equation holds at the origin by smoothness.

\medskip
\noindent\textbf{Step 2. Uniform estimates for the solution and the data.}\par\nobreak

Antisymmetry of $\Omega_N$ follows from that of $\Omega_T$.
By $0\le\chi\le1$ and \eqref{eq:potential-energy},
\[
\|\Omega_N\|_{L^2(B_1)}
\le\|\Omega_T\|_{L^2(B_1)}<\varepsilon,
\qquad
\|\nabla u_N\|_{L^2(B_1)}^2
=2\pi\int_T^\infty\frac{\chi(t/N)^2}{t^2}\,dt
\le\frac{2\pi}{T}.
\]
To estimate $U_N$, we use integration by parts. Since $\chi'$ is supported where $N\le t\le2N$,
\[
\left|\left(\chi(t/N)t^{-1/3}\right)'\right|
=\left|\frac{\chi'(t/N)}N t^{-1/3}
-\frac13\chi(t/N)t^{-4/3}\right|
\le Ct^{-4/3}.
\]
It follows that $|U_N(t)|\le Ct^{-1/3}$ for $t\ge T$, with a
constant independent of $N$. Hence
\[
\|u_N\|_{L^\infty(B_1)}\le CT^{-1/3},
\]
which, together with the gradient estimate, gives the required
uniform $W^{1,2}$ bound.

It remains to estimate the Hardy norm of the forcing. Differentiating
\eqref{eq:bn}, we obtain
\[
b_N'(t)
=-\frac{2\chi(t/N)}{t^3}
+\frac{2\chi'(t/N)}{Nt^2}
-\frac{\chi''(t/N)}{N^2t}.
\]
On the support of the cutoff derivatives, $N\le t\le2N$.
Thus $N^{-1}\le2t^{-1}$ there, and
\[
|b_N(t)|\le Ct^{-2},\qquad
|b_N'(t)|\le Ct^{-3},
\]
with a constant independent of $N$. By \eqref{eq:uniform-tail},
\[
\left|\int_t^\infty b_N(s)E(s)\,ds\right|
\le Ct^{-4/3}.
\]
Applying Lemma~\ref{lem:annular} to $F(t)=b_N(t)E(t)$, we conclude
that
\[
\|f_N^0\|_{h^1(\R^2)}\le CT^{-1/3}
\]
uniformly in $N$.

\medskip
\noindent\textbf{Step 3. Divergence of the Hessian norm.}\par\nobreak

For $t\le N$, we have $V_N(t)=t^{-1}E(t)$.
The tangential Hessian calculation in Step~4 of the proof of
Theorem~\ref{thm:main} therefore gives
\[
|D^2u_N(x)|\ge\frac1{r^2t}
\qquad\text{when }T+\log2\le t\le N.
\]
This interval corresponds to the annulus
$e^{T-N}\le r\le1/2$. Consequently,
\begin{equation*}\label{eq:smoothblowup}
\begin{aligned}
\|D^2u_N\|_{L^1(B_{1/2})}
&\ge2\pi\int_{T+\log2}^N\frac{dt}{t}\\
&=2\pi\log\left(\frac{N}{T+\log2}\right)
\to\infty.
\end{aligned}
\end{equation*}
Combining this with the uniform bounds in Step~2 completes the proof.
\end{proof}

\subsection*{Acknowledgements}\small
All authors are supported by the Young Scientist Program of the Ministry of Science and Technology
of China (No. 2021YFA1002200), the National Science Foundation of China (No. 12671095), the Taishan
Scholar Project and the Jiangsu Provincial Scientific Research Center of Applied Mathematics under Grant
No. 7707014040A. 

The authors would like to thank Prof.~Chang-Yu Guo and Prof.~Chang-Lin Xiang for introducing this problem to them and for valuable comments on the preliminary version of this draft. In particular, their suggestion on studying truncation of the construction in Theorem 1.1 leads to the formulation of Theorem 1.2 and Corollary 1.3.

\end{document}